\documentclass{article}

\usepackage{arxiv}

\usepackage[utf8]{inputenc} % allow utf-8 input
\usepackage[T1]{fontenc}    % use 8-bit T1 fonts
\usepackage{hyperref}       % hyperlinks
\usepackage{url}            % simple URL typesetting
\usepackage{booktabs}       % professional-quality tables
\usepackage{amsfonts}       % blackboard math symbols
\usepackage{nicefrac}       % compact symbols for 1/2, etc.
\usepackage{microtype}      % microtypography
\usepackage{lipsum}		% Can be removed after putting your text content
\usepackage{graphicx}

\usepackage{moreverb}

\usepackage{graphicx, subfigure}

\newcommand\BibTeX{{\rmfamily B\kern-.05em \textsc{i\kern-.025em b}\kern-.08em
T\kern-.1667em\lower.7ex\hbox{E}\kern-.125emX}}

\usepackage{graphicx}
\RequirePackage{fix-cm}
\usepackage{hyperref}
 \usepackage{epsfig} % for postscript graphics files
 \usepackage{amssymb}
 \usepackage{amsmath}
\usepackage{multirow}
\usepackage{epsfig} % for postscript graphics files
\input{epsf}
\usepackage{graphicx}
\usepackage{graphicx, subfigure}
\usepackage[utf8]{inputenc}
\usepackage[english]{babel}

\def\beq{\begin{equation}}
\def\eeq{\end{equation}}
\def\baq{\begin{eqnarray}}
\def\eaq{\end{eqnarray}}
\def\bal{\begin{align} }
\def\eal{\end{align} }
\def\bc{\begin{center}}
\def\ec{\end{center}}
\def\ds{\displaystyle}
\def\fai{\varphi}

\def\fai{\varphi}

\def\gr{\gamma_{r}}
\def\gr1{\gamma_{r1}}

\def\fai{\varphi}

\title{Kaczmarz Projection Algorithms in Moving Window : Performance Improvement via Extended Orthogonality \& Forgetting  \\  }
\author{ Alexander Stotsky \\
Department of Electrical Engineering \\
     Chalmers University of Technology  \\
     Gothenburg  SE - 412 96, Sweden  \\
	\texttt{alexander.stotsky@chalmers.se} \\
        \texttt{alexander.stotsky2@telia.com}  }
\hypersetup{
pdftitle={A template for the arxiv style},
pdfsubject={q-bio.NC, q-bio.QM},
pdfauthor={David S.~Hippocampus, Elias D.~Striatum},
pdfkeywords={First keyword, Second keyword, More},
}
\date{}
\begin{document}
\maketitle

\begin{abstract}
~New Kaczmarz algorithms with rank two gain update, extended orthogonality property and forgetting mechanism which
includes both exponential and instantaneous forgetting (implemented via a proper choice of the forgetting factor and the window size) are introduced and associated in this  report with well-known Kaczmarz algorithms with rank one update.
\end{abstract}


\maketitle

\section{Previous \& Related Work}
\label{introd} \vspace{-4pt}
\noindent
Kaczmarz projection method, \cite{astr}, \cite{ceg}   is robust and computationally efficient
algorithm (alternative to the least squares approach)
for estimation of the frequency contents of the oscillating signals
in a number of signal processing and system identification applications.
In the Kaczmarz algorithm  $ \ds \theta_k  = \theta_{k-1} -  \frac{\fai_k}{\fai^T_k ~ \fai_k}  ~ [ \fai^T_k~ \theta_{k-1} - y_k ]$  the output of the model $\fai^T_k~ \theta_{k}$  matches exactly the measured
signal $y_k = \fai^T_k~ \theta_{k}$  in each discrete step $k$, which can be obtained by multiplying  the parameter vector $\theta_k$ by the regressor $\fai^T_k$.  The algorithm has the remarkable orthogonality property
 $\fai_k^T ~ \tilde{\theta}_k = 0 $, where the vector of the parameter mismatch $\tilde{\theta}_k = \theta_k - \theta_*$ (difference between adjustable and true parameters $\theta_*$) is orthogonal to the regressor vector.
The square of the length of the harmonic regressor (which consists of the trigonometric
functions at different frequencies),  $\fai^T_k ~ \fai_k$ is constant and depends on the number of frequencies only.
The Kaczmarz algorithms were
successfully tested in different estimation problems in automotive applications, \cite{c}.
\\ Aiming for convergence rate improvement  of the Kaczmarz projection method
 the gain update algorithm was introduced in \cite{ave}.
 Further development was associated with introduction of the forgetting factor which
 has significant impact on the estimation performance, discounts exponentially  old measurements in the gain update and creates a virtual  moving window, \cite{st2012}.
 \\ The choice of the forgetting factor is associated with the trade-off between rapidity and accuracy of estimation.
 The gain matrix with forgetting factor  is associated in this paper with the
 window which is moving in time for  further improvements.  Introduction of the forgetting factor in the sliding window, \cite{liu},\cite{di} creates extended forgetting mechanism that includes both exponential and instantaneous forgetting and provides new opportunities for achievement of the trade-off between rapidity and accuracy.
 \\ The movement of the window is associated with data updating and downdating that results in recursive updates of the information matrix, which can occur sequentially, \cite{c} , \cite{zhang} or simultaneously, \cite{choi}.
Sequential updating and downdating results in computationally complex algorithms.
The gain update in known Kaczmarz algorithms, \cite{ave}, \cite{st2012} is associated with recursive rank one update, (similar to recursive least squares algorithms,  \cite{fom}, \cite{lju1})  whereas simultaneous updating/downdating is associated with  computationally efficient rank two update, \cite{sto2023}.
\\ New gain update which is explicitly associated with the movement of the sliding window of the size $w$ with
exponential forgetting is introduced in Kaczmarz algorithms. The parameter update law extends
the property of orthogonality where the output of the model matches measured data
in two end points of the weighted  window. The development is performed and systematically associated with well-known Kaczmarz algorithms with rank one update, \cite{st2012}.

\vspace{-10pt}
\section{Kaczmarz Algorithms in the Moving Window with Exponential Forgetting}
\label{mw} \vspace{-4pt}
\noindent
Suppose that the measured oscillating signal is
presented in the following form ${y}_k = \fai_k^{T} \theta_* $, where
$\fai_k$  is the harmonic regressor,  $\fai_k^{T} =  [ cos( q_0 k ) ~ sin(q_0 k ) ~ ...~ cos( q_h k ) ~ sin( q_h k ) ]$, $q_0,...q_h$ are the frequencies  and  $\theta_* $ is
the vector of unknown parameters, $k = 1,2,...$.
\\
The window of the size $w$ which is moving in time can be presented as rank two update of the information matrix (which is defined as the sum of the outer products of the regressor vectors)
where new observation is added (updating) and old observation is deleted (downdating).
The estimation performance is highly influenced by the forgetting factor, $ 0 < \lambda \le 1$ which discounts exponentially  old measurements and creates a virtual  window inside of the moving window.
In other words, the new data $\fai_{k}$ (with the largest forgetting factor which is equal to one) enter the window and the data with the lowest priority $ \tilde{\fai}_{k-w} = \sqrt{\lambda^w}~\fai_{k-w}$ leave the window in step $k$ as follows, \cite{sto2023}:
\begin{align}
 A_k & = \sum_{j = k -(w-1)}^{j=k} \lambda^{k-j} \fai_{j}~ \fai_{j}^{T} =  \lambda~A_{k-1} + Q_k~D~Q^{T}_k \label{r2update}
\end{align}
where $Q_k = [\fai_{k} ~ \tilde{\fai}_{k-w}]$,  $ D = diag[ 1, -1] $, $k \ge w + 1$.
\\
The gain matrix (as inverse of  $A_k$) which is derived by application of the matrix inversion lemma, \cite{hager} to the identity (\ref{r2update}) and the parameter update law can be written in the following form, \cite{sto2023},
\cite{jsp}:

\begin{align}
 & \Gamma_{k} = \frac{1}{\lambda} ~ [~ \Gamma_{k-1} - \Gamma_{k-1} ~ Q_k ~ S^{-1} ~ Q^T_k ~ \Gamma_{k-1}~]   \label{gk2} \\
  & \theta_k  = \theta_{k-1} - \Gamma_{k-1}  ~ Q_k ~ [Q^T_k ~ \Gamma_{k-1} ~ Q_k]^{-1} ~ [ Q^T_k~ \theta_{k-1} - \tilde{y}_k ] \label{ts1}
\end{align}
where $ S = \lambda~D + Q^T_k ~ \Gamma_{k-1} ~ Q_k $, $\tilde{y}^T_k = [y_k~ \sqrt{\lambda^w}y_{k-w}]$  is the synthetic output, provided that the matrix $ Q^T_k ~ \Gamma_{k-1} ~ Q_k $ is invertible. The vector of adjustable parameters $\theta_k$ in (\ref{ts1}) is constructed so that the output of the model matches exactly the measured signal  $Q^T_k~ \theta_{k} = \tilde{y}_k$ (see Section~\ref{introd}) and the following extended  orthogonality property $Q^T_k~\tilde{\theta}_k = 0 $ is valid for the parameter mismatch $\tilde{\theta}_k = \theta_{k} - \theta_*$.
\vspace{-4pt}
\section{Kaczmarz Algorithms with Rank One Update}
\label{kacu1} \vspace{-4pt}
\noindent
 Introduction of the forgetting factor allows to establish relationship between Kaczmarz
 algorithms with  rank two and rank one updates.
 Notice that $Q_k$ gets the following form
$Q_k = [\fai_{k} ~ 0]$ without downdating when $\lambda^w \to 0$. Straightforward   substitution of $Q_k$ defined above
in (\ref{gk2}) yields to following well-known Kaczmarz  algorithms, \cite{ave},\cite{st2012} :
\begin{align}
 \Gamma_{k} &= \frac{1}{\lambda} ~ [~ \Gamma_{k-1} - \frac{\Gamma_{k-1} ~ \fai_k ~ \fai^T_k ~ \Gamma_{k-1}~}{
 \lambda +  \fai^T_k ~ \Gamma_{k-1} \fai_k}~]   \label{gks} \\
 \theta_k &= \theta_{k-1} - \frac{\Gamma_{k-1} ~ \fai_k }{
 \fai^T_k ~ \Gamma_{k-1} \fai_k}~~[\fai^T_k \theta_{k-1} -  y_k]   \label{t1}
\end{align}
with the orthogonality property $\fai^T_k (\theta_k - \theta_*) = 0$.
\\ \\
Detailed description and comparisons of the Kaczmarz
 algorithms with  rank two and rank one updates together with their convergence properties are presented  in \cite{jsp} :
\\ \\
Stotsky A., Kaczmarz Projection Algorithms with Rank Two Gain Update,
 J Sign Process Syst, 2024.
 \\
    \url{https://doi.org/10.1007/s11265-024-01915-w}
\\ \\
where the performance  of new algorithms is also examined in the problem of estimation of the
grid events in the presence of significant harmonic emissions.

\vspace{-8pt}
\section{Conclusion}
\label{conc} \vspace{-4pt}
\noindent
The gain matrix with rank two update which allows  exponential weighting of the data inside of the moving window,
prioritizes recent measurements and improves estimation performance for fast varying changes of the signal
formed the basis for new Kaczmarz algorithms. New Kaczmarz adjustment law has extended orthogonality property in the sliding window.
The enhancement of the performance was achieved via a  proper choice of two adjustable parameters (small window size and forgetting factor) associated with fast forgetting.
\\ The results are valid under the assumption associated with invertibility of the matrix $Q^T_k ~ \Gamma_{k-1} ~ Q_k$
which can be seen as drawback of proposed algorithm.
For the sake of robustness the Kaczmarz algorithms with rank two update can be modified. The modification
can be associated with modification of well-known Kaczmarz  algorithm, \cite{astr}.

\end{document}